\documentclass{article}
\usepackage{overpic}
\usepackage{xcolor}
\usepackage{graphicx}
\usepackage{amsmath}
\usepackage{amsfonts}
\usepackage{amssymb}
\usepackage{amsmath,bm}
\usepackage{bbm}
\usepackage{pifont}
\usepackage{color}
\usepackage{amsopn, amstext}
\usepackage{mathrsfs}
\usepackage{hyperref}
\usepackage{textcomp}
\usepackage{url}
\usepackage{wasysym}
\usepackage[all]{xy}

\newcommand{\blind}{0}

\definecolor{dpur}{RGB}{163,12,85}

\newtheorem{theorem}{Theorem}[section]
\newtheorem{lemma}[theorem]{Lemma}
\newtheorem{corollary}[theorem]{Corollary}
\newtheorem{proposition}[theorem]{Proposition}

\newtheorem{definition}[theorem]{Definition}

\newtheorem{remark}[theorem]{Remark}

\numberwithin{equation}{section}

\DeclareMathOperator{\rank}{rank}

\def\E{\mathcal{E}}

\def\R{\mathbb{R}}

\def\dd{\mathrm{d}}

\def\s{\sigma}

\def\ra{\rightarrow}

\def\leqs{\leqslant}
\def\geqs{\geqslant}

\begin{document}

\def\spacingset#1{\renewcommand{\baselinestretch}%
{#1}\small\normalsize} \spacingset{1}


\if0\blind
{
  \title{\bf Hadamard-L\'{e}vy Theorems For Maps Taking Values in a Finite Dimensional Space}
  \author{Yacine Chitour\footnote{Laboratoire des Signaux et Syst\`{e}mes (L2S),  Universit\'{e} Paris-Saclay, CNRS,
  CentraleSup\'{e}lec,
 91190 Gif-sur-Yvette, France (email: yacine.chitour@centralesupelec.fr)}, ~Zhengping Ji\footnote{Chair for Dynamics, Control, Machine
Learning and Numerics (Alexander von Humboldt-Professorship), Department
of Mathematics, Friedrich-Alexander-Universit\"{a}t Erlangen-N\"{u}rnberg, 91058 Erlangen, Germany (email: zhengping.ji@fau.de)}, ~Emmanuel Tr\'{e}lat\footnote{Laboratoire Jacques-Louis Lions (LJLL), Sorbonne Universit\'{e}, CNRS, Universit\'{e} Paris Cit\'{e}, Inria, 75005 Paris,
France (email: emmanuel.trelat@sorbonne-universite.fr)}}
  \date{}
  \maketitle
} \fi

\if1\blind
{
  \bigskip
  \bigskip
  \bigskip
  \begin{center}
    {\LARGE\bf Title}
\end{center}
  \medskip
} \fi

\bigskip

\begin{abstract}
We propose global surjectivity theorems of differentiable maps based on second order conditions. Using the homotopy continuation method, we demonstrate that, for a $C^2$ differentiable map from a Hilbert space to a finite-dimensional Euclidean space, when its second-order differential has uniform upper and lower bounds, it has a global path-lifting property in the presence of singularities. This is then applied to the nonlinear motion planning  problem, establishing in some cases the well-posedness of the continuation method despite critical values of the endpoint maps.
\end{abstract}

\spacingset{1.45}

\section{Introduction}\label{S1}

The question of how surjectivity and injectivity of differentiable maps are governed by their differentials is a classical problem in analysis. We first recall the following fundamental result in nonlinear analysis \cite{amb,gr}.

\begin{theorem}[Local Surjectivity Theorem]
    Let $X$, $Y$ be Banach spaces, and $F:X\ra Y$ be a $C^1$ map. For $x\in X$, if $\dd F$ is invertible at $x$, then there exists a neighborhood $V$ of $F(x)$, such that $V\subset F(X)$.
\end{theorem}
To provide a global character to the above theorem, one needs to impose stronger conditions on the invertibility of $\dd F$. The first result in that direction was originally proposed by Hadamard \cite{ha} in $\R^n$, and then generalized by L\'{e}vy \cite{le} and Plastock \cite{pl} to Banach spaces.
\begin{theorem}[Hadamard-L\'{e}vy]\label{hl}
    Let $X$, $Y$ be Banach spaces, and $F:X\ra Y$ be a $C^1$ map with $\dd F|_x$ invertible for all $x\in X$. If there exists a continuous non-decreasing map $\beta:\R^+\ra\R^+\backslash\{0\}$ such that
    \begin{align}\label{0.0}
        \int_0^{+\infty}\frac{\dd s}{\beta(s)}=+\infty,\quad \|\dd F^{-1}|_x\|\leqs\beta(\|x\|), ~\forall x\in X,
    \end{align}
    then F is a $C^1$ diffeomorphism between $X$, $Y$.
\end{theorem}

A natural question is: to what extent such global surjectivity results can be extended to maps whose differentials are not globally invertible? This paper will explore this problem by exploiting the so-called {\it homotopy continuation method} (HCM) which is essential for the classical proof of Theorem \ref{hl} (one may refer to \cite{dem}\cite{pl}). In what follows we present the idea of homotopy continuation by sketching the proof of Hadamard-L\'{e}vy theorem, in order to introduce the methodology of our work.

The HCM can be traced back to the works of Poincar\'{e} \cite{po}, Klein \cite{kl}, Leray and Schauder \cite{ls}, and has found extensive application in numerical analysis \cite{aw}. Roughly speaking, it consists in embedding an equation into a parameterized family of equations so that the solution of this equation can be continuously deformed from simpler instances. The HCM has been applied in various schemes such as polynomial equation systems \cite{lt}, stability analysis \cite{sun},  motion planning \cite{ch02,se,chch}, because of its reliability, efficiency and regularity of solutions \cite{ri}.

For instance, in the proof of the Hadamard-L\'{e}vy theorem, the continuation method consists in establishing the path-lifting property of a map, as defined next.
\begin{definition}\label{d1.1}
	Let $X$, $Y$ be Banach spaces. Let $\gamma:[0,1]\ra Y$ be a $C^1$ curve with $\gamma(0)\in F(X)$. We say a $C^1$ map $F:X\ra Y$ lifts $\gamma$ if there exists a $C^1$ curve $\varGamma:[0,1]\ra X$, such that $F\circ\varGamma(s)=\gamma(s)$ for all $s\in[0,1]$.
\end{definition}

Using the above definition, we next present the main idea of proving Theorem \ref{hl} via continuation, as preparation for its generalization. By definition \ref{d1.1}, if $Y$ is path-connected and $F$ lifts any curve on $Y$, then $F$ is surjective, since for any $y\in Y$, by the path-connectedness there exists a curve $\gamma$ with $\gamma(0)\in F(X)$ and $\gamma(1)=x$, while the lift of $\gamma$ gives $F(\varGamma(1))=y$, $\varGamma(1)\in X$. Hence it suffices to demonstrate the path-lifting property of $F$ for showing its surjectivity in Theorem \ref{hl}. 

Let $y\in Y$, and $\gamma$ be a $C^1$ curve  with $\gamma(0)\in F(X)$ and $\gamma(1)=y$. Consider the following Wazewski equation \cite{wa} of $\varGamma:[0,1]\ra X$ with respect to $F$ and $\gamma$
    \begin{align}\label{0.2}
        \dot{\varGamma}(s)=(\dd F|_{\varGamma(s)})^{-1}\dot{\gamma}(s)
    \end{align}
    with initial condition $F(\varGamma(0))=\gamma(0)$. It is shown in \cite{dem,pl} that under condition \eqref{0.0}, the flow of this equation is complete over $s\in[0,1]$, which means that the map $F$ lifts $\gamma$, $y=F(\varGamma(1))$, with $\varGamma(1)$ given by the terminal value of the dynamics. On the other hand, since $F$ is a local homeomorphism due to the invertibility of $\dd F$, the lift of $\gamma$ is unique. Therefore the map $F$ satisfying \eqref{0.0} is a diffeomorphism.

    Inspired by the above continuation method of proving the Hadamard-L\'{e}vy theorem, we try to extend it with the presence of singularities in the case when $X$ and $Y$ are not of the same dimensionality. We first rephrase the continuation under general nonlinear settings. Let $F:X\ra M$ be a differentiable map from  a Hilbert space $X$ to  a smooth Riemannian manifold $M$. Define the \emph{Gramian} of $F$ at $u\in X$ as
    \begin{align}\label{0.0.3}
        G(u):=\dd F|_u\dd F|_u^*,
    \end{align}
    then $G(u)$ is a linear map from $T_{F(u)}^*M$ to $T_{F(u)}M$. Using this notion, consider the following nonlinear partial differential equation  
    \begin{align}\label{1.1}
		\frac{\partial u(s)}{\partial s}=\dd F|^*_{u(s)}G(u(s))^{-1}\dot{\gamma}(s)
	\end{align}
	over $s\in[0,1]$ with initial condition $u(0)=u^0$. Then the surjectivity of the map $F$ can be established by studying the completeness of the above equation, as is shown next.
    
    The following lemma is a generalization of the results in \cite{ch06,su93}, in which the authors studied the surjectivity of a specific  function $F$, namely the endpoint map of nonholonomic control systems. 
\begin{lemma}\label{l1.0}
	Let $X$ be a Hilbert space and $(M,g)$ be a path-connected smooth Riemannian manifold. Let $F:X\ra M$ 
 be a $C^2$ differentiable map such that $\dd F|_u$ is surjective for all $u\in X$. If there exists  $C>0$ such that 
 \begin{align}\label{0.11}
     \|G(u)^{-1}\|\leqs{C}(1+{\|u\|_X})^2, \quad  \forall u\in X,
 \end{align}
 where $\|\cdot\|_X$ is the norm on $X$ and $\|\cdot\|$ is the operator norm of a linear map from $T_{F(u)}M$ to $T^*_{F(u)}M$, then for any $u^0\in X$ and any $C^1$ curve $\gamma:[0,1]\ra M$ satisfying $\gamma(0)=F(u^0)$, there exists a global solution to the equation \eqref{1.1}, and $F$ is surjective.
\end{lemma}

The proof consists in using the Gr\"{o}nwall's lemma to show the global existence of solutions to \eqref{1.1},  which can be easily constructed following the aforementioned references.

The equation \eqref{1.1} is called the \emph{path-lifting equation} (PLE). It is the least-square version of \eqref{0.2} in the proof of the Hadamard-L\'{e}vy theorem, ensuring the existence of lifts.  In the nonlinear case, when \eqref{0.11} holds,  $\rank(\dd F|_u)=\dim(M)$ for all $u\in X$. Then for any $s\in[0,1]$, \eqref{1.1} gives the least-square solution to the equality $\dd F|_{u(s)}(\frac{\partial  u(s)}{\partial s})=\dot{\gamma}(s)$, since $\dd F$ is not assumed to be injective.

On the other hand, when $\dd F$ is not everywhere invertible on $X$, the above path-lifting is not applicable in general as before. Adopting the definition of the Gramian as in \eqref{0.0.3}, we denote by $\tilde{S}$ the \emph{singular set} of $F$, which is the closed set defined as
\begin{align}\label{0.0.0}
 \tilde{S}:=\{u\in X\:|\:\det\big(G(u)\big)=0\}
\end{align}
and denote by $S:=F(\tilde{S})$ the set of \emph{critical values} of $F$. The global existence of the solutions of equation \eqref{1.1} is not guaranteed when $u(s)$ is approaching $\tilde{S}$, hence the applicability of the continuation method in Lemma \ref{l1.0} as well as in the Hadamard-L\'{e}vy theorem is questionable in the presence of singularities, which gives rise to the main subject of this paper.

Concerning the maps with non-invertible first order differentials, there exist results on their invertibility in finite dimension \cite{ava}, and in Banach spaces with conditions on the corank of the critical values \cite{am}. In the case where $F$ is the endpoint map corresponding to a control-affine system, \cite{sch} proposed a regularized path-lifting equation (RPLE) for studying its surjectivity in the motion planning problem in the presence of singularities, and \cite{ji} proved the convergence of its solutions under appropriate conditions.

In this paper, we  present Hadamard-L\'{e}vy type theorems which provide second-order sufficient conditions for the path-lifting property and surjectivity  of differentiable maps from a Hilbert space to a finite-dimensional Euclidean space, in the presence of singularities. The idea is to adopt the continuation method as in Lemma \ref{l1.0} and to investigate the global existence of the solution of the path-lifting equation \eqref{1.1} via spectral analysis of the Gramian matrix, estimating its least eigenvalue and the associated switching functions. This allows us to establish surjectivity and the path-lifting property even when the first-order differential fails to be invertible.

The rest of this paper is structured as follows. 
Section \ref{S0} presents the main theorems, and sketches the idea of the proof by showing the basic estimates on the solutions of the PLE which helps establishing the local well-posedness.
In Section \ref{S2} we construct the estimates on the derivatives of the spectral of the Gramian along the PLE, and in Section \ref{S3} we use these estimates to verify the global path-lifting property and to prove the main theorems. Then in Section \ref{S4} we apply these results to the motion planning problem of nonlinear control systems, showing that under proposed second-order conditions, the continuation method remains robust and well-posed, even when encountering critical values in the endpoint map.

\section{Main Results}\label{S0}

\subsection{Notations and main theorems}

Throughout, let $X$ be a Hilbert space and $\|\cdot\|_X$ be the norm induced by its inner product. Let $F:X\ra \R^n$ be a second-order differentiable map and $u\in X$, $z\in\R^n$, and consider the nonlinear map 
\begin{align}\label{0.1.3}
\begin{array}{rl}
     \phi_z:X\!\!\!&\ra X^*\\
u\!\!\!&\mapsto \dd F|_u^*z   
\end{array} 
\end{align}
called the \emph{switching function} of $F$ at $u$ with respect to $z$. By the standard isomorphism between $X$ and $X^*$, one can also view $\phi$ as taking values in $X$. Moreover, in the case when $u\in\tilde{S}$, there exists a nonzero $z\in\R^n$ such that $\phi_z(u)=0$. 

Further, for any $u\in X$ and $z\in\R^n$ we define a bilinear map 
\begin{align}\label{0.1.1}
\begin{array}{rl}
     z^*\dd^2 F|_u:X\times X&\ra ~\R\\
    (v,w)&\mapsto~\dd\phi_z|_u(v)(w),
\end{array}
\end{align}
where $\dd\phi_z$ is the differential of the switching function. 

For any $u\in X$, consider the Gramian matrix of $F$ at $u$ defined as in \eqref{0.0.3}. Denote the ordered eigenvalues of $G(u)$ as  $\lambda_1(u)\leqs\lambda_2(u)\leqs\cdots\leqs\lambda_n(u)$. Obviously, by definition \eqref{0.0.0}, $u\in\tilde{S}$ if and only if $\lambda_1(u)=0$.

In this paper we consider only corank-1 singularities and {we make the following assumption which will hold throughout, which guarantees the simplicity of the least eigenvalue of the Gramian.
\begin{enumerate}
    \item[\bf(A)] $\exists\lambda_0>0,~\text{s.t.}~\lambda_i(u)\geqs\lambda_0,~\forall i=2,\ldots,n, ~\forall u\in X.$
\end{enumerate}}

With the above assumption and definitions \eqref{0.1.3}\eqref{0.1.1}, the main results of this paper are the following theorems.

\begin{theorem}\label{t1}
	Let $F:X\ra\R^n$ be a $C^2$ map satisfying the following assumptions.
    \begin{itemize}
    \item $F$ has no singularities and Assumption (A) holds everywhere;
    \item there exist positive constants $C$, $R$, and a continuous non-decreasing map $\xi:\R^+\ra\R^+\backslash\{0\}$ such that for any $u\in X$ with $\|u\|_X\geqs R$ and any $z\in \R^n$ with $\|z\|=1$, 
\begin{subequations}\label{eq:parent}
	\begin{align}
		&|z^*\dd^2F|_u(\cdot,\cdot)|\leqs C;\label{eq:1}\\
		&\|\dd F|_u^*(z)\|\big|z^*\dd^2F|_u\big(\phi_z(u),\phi_z(u)\big)\big|\geqs \frac{\|\phi_z(u)\|_X^2}{\xi(\|u\|)^2},\quad\int_R^{+\infty}\!\!\frac{\dd s}{\xi(s)}=\infty\label{eq:2} 
	\end{align}    
\end{subequations}
with $\phi_z(u)$ defined as in \eqref{0.1.3}.
\end{itemize}
Then the map $F$ is surjective. 
\end{theorem}

\begin{remark}\label{r2.0.0}
    One special case of the condition \eqref{eq:2} is that when there exist constants $K_1$, $K_2\geqs0$ and $\alpha\in(0,1)$ such that 
    \begin{align}\label{0.0.01}
\big|z^*\dd^2F|_u\big(\phi_z(u),\phi_z(u)\big)\big|\geqs K_1\frac{\|\phi_z(u)\|_X^2}{\|u\|^{1-\alpha}_X};\quad{\|\dd F|_u^*(z)\|_X}\geqs \frac{K_2}{\|u\|_X^{1+\alpha}}
    \end{align}
    holds for all $u\in X$ with $\|u\|_X\geqs R$. In this case $\xi(s)=s$, and one sees that by taking the function $\beta(s)=s^{1+\alpha}$,  which does not satisfy the condition verified by $\beta$ in \eqref{0.0}, Theorem \ref{t1} becomes an extension of the original Hadamard-L\'{e}vy theorem (for maps taking value in Euclidean spaces) at the price of extra assumptions.
\end{remark}

On the other hand, since the first-order differential of the map $F$ is supposed to be full rank in \eqref{eq:2}, the above theorem does not concern the size of the preimage of the critical value of $F$. We further propose the following theorem on the path-lifting property of $F$, which allows singularity at the terminal of the lifted path.

\begin{theorem}\label{t1.1}
Let $F:X\ra\R^n$ be a $C^2$ map and and let $\tilde{S}$ be its 
singular set defined as in \eqref{0.0.0}. Assume that $F$ satisfies the following.
 \begin{itemize}
    \item Assumption (A) holds everywhere;
    \item there exist positive constants $C$, $R$, such that for any $u\in X\backslash\tilde{S}$ with $\|u\|_X\geqs R$ and any $z\in \R^n$ with $\|z\|=1$, 
\begin{align}\label{eq1:2} 
		|z^*\dd^2F|_u(\cdot,\cdot)|\leqs C;\quad |z^*\dd^2F|_u\big(\phi_z(u),\phi_z(u)\big)|\geqs K\|\phi_z(u)\|_X^2
	\end{align} 
with $\phi_z(u)$ defined as in \eqref{0.1.3}.
\end{itemize}
Then for any $C^2$ curve $\gamma:[0,1]\ra \R^n$ satisfying $\gamma(s)\notin S$ for all $s\in[0,1)$, $F$ lifts $\gamma$. 
\end{theorem}

We will prove the above theorems by investigating the global existence of the solution of the  PLE \eqref{1.1} with $u^0\notin \tilde{S}$, when the curve $\gamma$ to be lifted satisfies $\gamma(s)\notin S$ for $s\in [0,1)$. As explained in Section \ref{S1}, if the PLE is globally well-posed, then the path-lifting property implies the surjectivity of the map $F$.

In the rest of this section, we will introduce the notations needed to prove the above main results, then establish the local well-posedness of the PLE and present the key estimate for showing global existence of the solution, sketching the main idea of the proofs of Theorems \ref{t1} and \ref{t1.1}.


For any $s\in[0,1)$ such that $u(s)$ exists,  the Gramian matrix (defined as in \eqref{0.0.3}) at $u(s)$ and its ordered (non-negative) eigenvalues are denoted by 
\begin{align}\label{1.0.1}
   G(s):=G(u(s))=\dd F|_{u(s)}\dd F|_{u(s)}^*\in\R^{n\times n},
\end{align}
and
\begin{align}\label{1.0.2}
    \mathrm{Spec}(G(s))=\{\lambda_1(s),\ldots, \lambda_n(s)\}, \quad \lambda_1(s)\leqs\cdots\leqs\lambda_n(s)
\end{align}
Let $z_1(s),\ldots,z_n(s)$ be eigenvectors of norm one, i.e., 
\begin{align}\label{1.0.3}
    G(s)z_i(s)=\lambda_i(s)z_i(s), \quad \|z_i(s)\|_{\R^n}=1, \quad i=1,\ldots,n.
\end{align}
Since $s\mapsto u(s)$ is continuous where it is defined, and since $u\mapsto G(u)$ is $C^1$ continuous with respect to $u\in X$ due to the $C^2$ continuity of $F$, we see that the real value functions 
$s\mapsto \lambda_i(\cdot)$ $i=1,\ldots,n$ are continuous and one can choose vector-valued maps $s\mapsto z_i(\cdot)$, $i=1,\ldots,n$ to be also continuous, as long as $u(\cdot)$ is defined. In particular, as Assumption (A) guarantees the simplicity of $\lambda_1(s)$, we have the uniqueness of $z_1(s)$ up to a sign.

By definition of the Gramian and the switching function \eqref{0.1.3}, for each unit eigenvector along the solution of the PLE, we have
\begin{align}\label{0.0.1}
    \forall s\in[0,1),~ \|\phi_{z_i(s)}\big(u(s)\big)\|=\|\dd F|_{u(s)}^*z_i(s)\|_X=\sqrt{\lambda_i(s)}, \quad i=1,\ldots,n.
\end{align}
When $u(s)$ exists and $u(s)\notin\tilde{S}$, we define the \emph{normalized switching functions along $u(s)$} as  follows:
\begin{align}\label{0.0.2}
    v_i(s):=\frac{\dd F|^*_{u(s)}z_i(s)}{\sqrt{\lambda_i(s)}},\quad i=1,\ldots,n.
\end{align}
By \eqref{0.0.1}, $v_i(s)\in X$, $i=1,\ldots,n$ are of unit length for any $s$ such that $u(s)$ exists outside of the singular set. 


\subsection{Brief outlines of the proofs}

With the above notations, we next outline the structure of the proofs of Theorems \ref{t1} and \ref{t1.1}. First, notice that if the initial condition $u^0$ does not belong to $\tilde{S}$ then one has local existence of the solution of the PLE \eqref{1.1} by the Cauchy-Lipschitz theorem. Indeed since $F$ is $C^2$ continuous with respect to $u$, then the right-hand side of 
\eqref{1.1} is $C^1$ continuous with respect to $u\in X\setminus \tilde{S}$. 

 With the notations in \eqref{1.0.3}, we start by writing
	\begin{align}\label{3}
		\dot{\gamma}(s)=\sum_{i=1}^na_i(s)z_i(s),
	\end{align}
    for $s$ belonging to the interval of existence of  solutions of \eqref{1.1}, denoted by $I$. Due to the continuity of $z_i(s)$ we obtain the continuous functions $a_i(s)$ over $s\in I$, $i=1,\ldots,n$. By the $C^2$ continuity of $\gamma(s)$ over $s\in[0,1]$, $\big\{|a_i(s)|~\big|~i=1,\ldots,n,~s\in I\big\}$ is bounded (independently of $I$). 
    Decomposing the right hand side of the PLE \eqref{1.1} with respect to the unit eigenvectors of the Gramian as
    \begin{align}\label{3.1.1}
        \frac{\partial u(s)}{\partial s}=\sum_{i=1}^n\frac{a_i(s)\dd F|_{u(s)}^*z_i(s)}{\lambda_i(s)}=\sum_{i=1}^n\frac{a_i(s)}{\sqrt{\lambda_i(s)}}v_i(s),
    \end{align}
	with $v_i(s)$, $i=1,\ldots,n$ defined as in \eqref{0.0.2}, we have the following lemma.
	\begin{lemma}\label{l3.2}
		Let $u(\cdot)$ be a solution to the PLE \eqref{1.1} with initial condition $u(0)=u^0\notin \tilde{S}$ defined on some open neighborhood $I$ of $0$. Then, there exists a constant $\bar{C}>0$ (only depending on $\gamma$), such that
		\begin{align}\label{3.00}
			\Big\|\frac{\partial u(s)}{\partial s}\Big\|_X\leqs \frac{|a_1(s)|}{\sqrt{\lambda_1(s)}}+\bar{C},
		\end{align}
		where $\lambda_1(s)$ is defined as in \eqref{1.0.2} and $a_1(s)$ as in \eqref{3}.

	\end{lemma}
	{\it Proof.} 
Since $z_i(s)$ is the eigenvector of $G(s)$ corresponding to $\lambda_i(s)$, by \eqref{3.1.1} and \eqref{0.0.1} we have
	\begin{align}\label{3.111}
		\Big\|\frac{\partial u(s)}{\partial s}\Big\|_X=\Big\|\sum_{i=1}^n\frac{a_i(s)}{\sqrt{\lambda_i(s)}}v_i(s)\Big\|_X\leqs \frac{|a_1(s)|}{\sqrt{\lambda_1(s)}}+ \sum_{i=2}^n\frac{|a_i(s)|}{\sqrt{\lambda_i(s)}}
	\end{align}
	 On the other hand, since (A) holds, for any $s\in[0,1)$ such that $u(s)$ exists,
	$$
	\sum_{i=2}^n\frac{|a_i(s)|}{\sqrt{\lambda_i(s)}}\leqs
 \frac{(n-2)\|\dot \gamma(s)\|_{\R^n}}{\sqrt{\lambda_0}},
	$$
	which proves the lemma by taking 
 $ \bar{C}:= \frac{(n-2)\max\limits_{s\in [0,1]}\big\|\dot \gamma(s)\big\|_{\R^n}}{\sqrt{\lambda_0}}$.
	\hfill{$\Box$}

In the sequel, we derive estimates of $\sqrt{\lambda_1(s)}$ and $a_1(s)$ defined as in \eqref{1.0.2}\eqref{3}, and compute their ratio to establish a uniform bound for the solution of the PLE. To be specific, we will prove the main theorems with the following procedure:
\begin{enumerate}
    \item show that both $\lambda_1(s)$ and $a_1(s)$ are differentiable outside of singularities (Proposition \ref{p3.1}, \ref{p3.2});
    \item build coupled differential equations of $a_1(s)$ and $\sqrt{\lambda_1(s)}$, show the integrability of $\Big|\frac{a_1(s)}{\sqrt{\lambda_1(s)}}\Big|$ outside of singularities  in order to prove Theorem \ref{t1};
    \item establish a differential equation of $g(s):=\frac{a_1(s)}{\sqrt{\lambda_1(s)}}$ to show its absolute integrability on singularities, proving Theorem \ref{t1.1}.
\end{enumerate}

\section{Differentiation and Estimates of the Spectrum}\label{S2}

With the notations introduced in Section \ref{S0}, in this section we estimate the differentials of $\sqrt{\lambda_1(s)}$ and $a_1(s)$ when $\lambda_1(s)>0$, in preparation for proving the main results on global well-posedness of the PLE following Lemma \ref{l3.2}.
 Define 
    \begin{align}\label{3.0.0}
        \s_0:=\sup\big\{\sigma\in[0,1)~\big|~\lambda_1(s)>0, ~\forall s\in[0,\sigma]\:\big\},\quad I:=[0,\sigma_0).
    \end{align}
If $\varliminf\limits_{s\ra\sigma_0}\lambda_1(s)>0$, then by Lemma \ref{l3.2}, the PLE \eqref{1.1} is globally well-posed and there is nothing to prove. Therefore, from now on, we assume that 
    \begin{enumerate}
        \item[\bf(A')] $\varliminf\limits_{s\ra\sigma_0}\lambda_1(s)=0$ and $\lambda_1(s)>0$ for $s\in I$. 
    \end{enumerate}

As discussed in the previous section, the eigenvalues $\lambda_i$, $i=1,\ldots,n$ are continuous functions of $s$; moreover, since $\lambda_1$ is simple, then $s\mapsto \lambda_1(s)$ is differentiable on  $I$ by the Danskin's theorem, cf. \cite[Appendix B, Section B.5]{ber}. As a consequence, further $s\mapsto z_1(s)$ is also differentiable on $I$.

The differential of the Gramian at $u\in X$ acting on two vectors $z_1$, $z_2\in \R^n$ is equal to
    \begin{align}\label{4.0.2}
        \begin{array}{rl}
             z_1^*\dd G|_uz_2:X\!\!&\ra ~\R\\
        v\!\!&\mapsto ~z_1^*\dd^2F|_u(\dd F|_u^*(z_2),v)+z_2^*\dd^2F|_u(\dd F|_u^*(z_1),v).
        \end{array}
    \end{align}
    with $z_i^*\dd^2 F|_u(\cdot,\cdot)$, $i=1,2$ defined as in \eqref{0.1.1}.
    
	Based on this, we make the following estimate on the derivative of the least eigenvalue of the Gramian along the  solution of the PLE.

	\begin{proposition}\label{p3.1}
 
 There exists $C_0\geqs0$ such that, for every $s\in I$, $\Big|\frac{\dd\lambda_1}{\dd s}(s)\Big|\leqs C_0$  and
    \begin{equation}\label{3.5.0}
        \frac{\dd \sqrt\lambda_1}{\dd s}(s)=\frac{a_1(s)}{\sqrt{\lambda_1(s)}}h(s)+f(s),
     \end{equation}
where $a_1(s)$ is defined as in \eqref{3}, and
    \begin{subequations}\label{eq2:parent}
        \begin{align}
			&h(s):={z_1(s)^*}\dd^2F|_{u(s)}\Big(v_1(s),v_1(s)\Big),\label{eq2:1}\\
            &f(s):=\sum_{i=2}^n\frac{a_i(s)}{\sqrt{\lambda_i(s)}}z_1(s)^*\dd^2F|_{u(s)}(v_i(s),v_1(s)),\label{eq2:2}
		\end{align}
    \end{subequations}
	with $z_i(s)$ defined as in \eqref{1.0.3} and $v_i(s)$ as in \eqref{0.0.2}. Moreover, $f$ and $\lambda_1$ are bounded over $I$.  
 
	\end{proposition}
	{\it Proof.}
	Let $u(\cdot)$ be the solution of \eqref{1.1}. 
	Then, by the decomposition \eqref{3.1.1} and the expression \eqref{4.0.2}, differentiating the least eigenvalue of the Gramian yields
for $s\in I$ that 
\begin{align}
		\frac{\dd \lambda_1}{\dd s}(s)=&\frac{\dd z_1(s)^*G(s)z_1(s)}{\dd s}\nonumber\\
		=&z_1(s)^*\frac{\dd G(s)}{\dd s}z_1(s)\nonumber\\
		=&z_1(s)^*\dd G_{u}|_{u(s)}\Big(\frac{\partial u(s)}{\partial s}\Big)z_1(s)\nonumber\\
        =&2z_1(s)^*\dd^2F|_{u(s)}\Big(\frac{\partial u(s)}{\partial s}, \dd F|_{u(s)}^*z_1(s)\Big)\nonumber\\
        =&{2}{\sqrt{\lambda_1(s)}}z_1(s)^*\dd^2 F|_{u(s)}\Big(\sum_{i=1}^n\frac{a_i(s)}{\sqrt{\lambda_i(s)}}v_i(s),v_1(s)\Big)\nonumber\\
        =&2a_1(s)z_1^*(s)\dd^2F|_{u(s)}(v_1(s),v_1(s))+\sum_{i=2}^n2a_i(s)\sqrt{\frac{\lambda_1(s)}{\lambda_i(s)}}z_1(s)^*\dd^2F|_{u(s)}(v_i(s),v_1(s))\nonumber\\
		=&\:2a_1(s)h(s)+2f(s)\sqrt{\lambda_1(s)},\label{3.1}
	\end{align}
	where $h(s)$ and $f(s)$ are defined as in \eqref{eq2:parent}. As $\frac{\dd\sqrt{\lambda_1}}{\dd s}(s)=\frac{1}{2\sqrt{\lambda_1(s)}}\frac{\dd\lambda_1}{\dd s}(s)$, \eqref{3.5.0} is established.
    
    By Condition \eqref{eq:1}, $|h(s)|\leqs C$ for $s\in I$ and by the $C^2$-continuity of $\gamma(s)$, $|a_1(s)|\leqs\|\dot\gamma(s)\|_{\infty}$ is bounded for $s\in I$. Using \eqref{eq2:2}, one has on $I$ that
    \begin{align*}
       |f(s)|\sqrt{\lambda_1(s)}
        \leqs2C\sum_{i=2}^n|a_i(s)|\sqrt{\frac{\lambda_1(s)}{\lambda_i(s)}}\leqs&\sup\limits_{s\in[0,\sigma_0),~i=2,\cdots,n}2(n-1)C|a_i(s)|
    \end{align*}
    where  the inequality is due to \eqref{eq:1} and the fact that $\lambda_1(s)\leqs\lambda_i(s)$ on $I$, $i=2,\ldots,n$. Therefore, $|f(s)|\sqrt{\lambda_1(s)}$ is bounded on $I$. Together with the uniform boundedness of $a_1(s)$ and $h(s)$, this implies by the estimate \eqref{3.1} that $|\frac{\dd \lambda_1}{\dd s}(s)|$ is bounded over $I$. Since $I$ is bounded (subset of $[0,1]$), it follows that $\lambda_1$ is bounded over $I$.

   Finally, we show the boundedness of $f(s)$. Since  $\lambda_i(s)\geqs\lambda_0$ for all $i=2,\ldots,n$ and $s\in I$ according to assumption (A) and \eqref{eq:1}, we have
    \begin{align*}
        |f(s)|\leqs\sum_{i=2}^n\sup\limits_{s\in[0,\sigma_0)}\frac{|a_i(s)|}{\sqrt{\lambda_i(s)}}C,
    \end{align*}
    and the conclusion follows since the $a_i$'s are bounded and by using assumption (A).
	\hfill{$\Box$}




\begin{proposition}\label{p3.2}
	The coefficient $a_1$ defined as in \eqref{3} verifies the following o.d. on $I$, 
	\begin{equation}\label{4.2}
		\frac{\dd a_1}{\dd s}(s)= \frac{a_1(s)}{\sqrt{\lambda_1(s)}}f(s)+O(1),
	\end{equation}
    where $O(1)$ denotes a function bounded over $I$. 
\end{proposition}
{\it Proof.} 
For $s\in I$, recall that  $G(s)z_1(s)=\lambda_1(s)z_1(s)$ and since  both $\lambda_1(s)$ and $z_1(s)$ are differentiable on $I$, differentiating this equality yields
\begin{align}\label{4.0}
	\frac{\dd z_1}{\dd s}(s)=&\Big(G(s)-\lambda_1(s)I_n\Big)\Big|_{V_1^{\perp}}^{-1}\Big(\frac{\dd G(s)}{\dd s}-\frac{\dd \lambda_1}{\dd s}(s)I_n\Big)z_1(s)\nonumber\\
	=&\Big(G(s)-\lambda_1(s)I_n\Big)\Big|_{V_1^{\perp}}^{-1}\left(\dd G|_{u(s)}\Big(\frac{\partial u}{\partial s}(s)\Big)z_1(s)-\frac{\dd \lambda_1}{\dd s}(s)z_1(s)\right)\nonumber\\
	=&\Big(G(s)-\lambda_1(s)I_n\Big)\Big|_{V_1^{\perp}}^{-1}\dd G|_{u(s)}\Big(\frac{\partial u}{\partial s}(s)\Big)z_1(s)
\end{align}
where the map $\Big(G(s)-\lambda_1(s)I_n\Big)\Big|_{V_1^{\perp}}^{-1}$ is defined as
\begin{align}\label{4.02}
    \Big(G(s)-\lambda_1(s)I_n\Big)\Big|_{V_1^{\perp}}^{-1}:=\sum_{i=2}^n\frac{1}{\lambda_i(s)}P_i(s),
\end{align}
with $P_i(s)$ the canonical projector of $G(s)$ with respect to $\lambda_i(s)$, $i=2,\ldots,n$. 
The last equality in \eqref{4.0} is due to the fact that $\Big(G(s)-\lambda_1(s)I_n\Big)\Big|^{-1}_{V_1^{\perp}}z_1(s)=0$.

Expressing $\frac{\partial u}{\partial s}$ by the PLE and adopting the expression of $\dd G$ in \eqref{4.0.2}, by \eqref{4.0} we have
\begin{align}\label{4.01}
    \Big\langle \zeta,\frac{\dd z}{\dd s}(s)\Big\rangle=&\Big\langle \zeta,\Big(G(s)-\lambda_1(s)I_n\Big)\Big|_{V_1^{\perp}}^{-1}\dd G_u|_{u(s)}\Big(\sum_{j=1}^n\frac{a_j(s)}{\sqrt{\lambda_j(s)}}v_j(s)\Big)z_1(s)\Big\rangle\nonumber\\
    =&z_1^*(s)\dd^2F|_{u(s)}\Big(\dd F|^*_{u(s)}\Big(\Big(G(s)-\lambda_1(s)I_n\Big)\Big|_{V_1^{\perp}}^{-1}\zeta\Big),\sum_{j=1}^n\frac{a_j(s)}{\sqrt{\lambda_j(s)}}v_j(s)\Big)\nonumber\\
    ~&+\Big\langle\zeta,\Big(G(s)-\lambda_1(s)I_n\Big)\Big|_{V_1^{\perp}}^{-1}\dd^2 F|_{u(s)}\Big(\dd F|^*_{u(s)}z_1(s),\sum_{j=1}^n\frac{a_j(s)}{\sqrt{\lambda_j(s)}}v_j(s)\Big)\Big\rangle
\end{align}
for any $\zeta\in\R^n$.

Then we use \eqref{4.01} to derive the function $a_1(s)$ in the decomposition \eqref{3}.
\begin{align}\label{4.1}
	\frac{\dd a_1}{\dd s}(s)
	=&\frac{\dd}{\dd s}\big\langle\dot{\gamma}(s),z_1(s)\big\rangle\nonumber\\
	=&\big\langle\ddot{\gamma}(s),z_1(s)\big\rangle+\left\langle\dot{\gamma}(s),\frac{\dd z_1}{\dd s}(s)\right\rangle\nonumber\\
	=&\big\langle\ddot{\gamma}(s),z_1(s)\big\rangle+\left\langle\sum_{i=1}^na_i(s)z_i(s),\Big(G(s)-\lambda_1(s)I_n\Big)\Big|_{V_1^{\perp}}^{-1}\Big(\dd G_u|_{u(s)}\big(\frac{\partial u}{\partial s}(s)\big)z_1(s)\Big)\right\rangle\nonumber\\
	=&\big\langle\ddot{\gamma}(s),z_1(s)\big\rangle+z_1^*(s)\dd^2F|_{u(s)}\Big(\dd F|^*_{u(s)}\big(G(s)-\lambda_1(s)I_n\big)\Big|_{V_1^{\perp}}^{-1}\Big(\sum_{i=2}^n a_i(s)z_i(s)\Big),\sum_{j=1}^n\frac{a_j(s)}{\sqrt{\lambda_j(s)}}v_j(s)\Big)\nonumber\\
	&+\left\langle\Big(G(s)-\lambda_1(s)I_n\Big)\Big|_{V_1^{\perp}}^{-1}\Big(\sum_{i=2}^n a_i(s)z_i(s)\Big),\dd^2F|_{u(s)}\Big(\dd F|^*_{u(s)}z_1(s),\sum_{j=1}^n\frac{a_j(s)}{\sqrt{\lambda_j(s)}}v_j(s)\Big)\right\rangle\nonumber\\
	=&\big\langle\ddot{\gamma}(s),z_1(s)\big\rangle+\frac{a_1(s)}{\sqrt{\lambda_1(s)}}\left(\sum_{i=2}^n\frac{a_i(s)\sqrt{\lambda_i(s)}}{\lambda_i(s)-\lambda_1(s)}z_1^*(s)\dd^2 F|_{u(s)}\big(v_i(s),v_1(s)\big)\right)\nonumber\\
	&+\sum_{i,j=2}^n\frac{a_i(s)a_j(s)\sqrt{\lambda_i(s)}}{\big(\lambda_i(s)-\lambda_1(s)\big)\sqrt{\lambda_j(s)}}z_1^*(s)\dd^2F|_{u(s)}\big(v_i(s),v_j(s)\big)\nonumber\\
	&+\sum_{i,j=2}^n\frac{a_i(s)a_j(s)\sqrt{\lambda_1(s)}}{\big(\lambda_i(s)-\lambda_1(s)\big)\sqrt{\lambda_j(s)}}z_i^*(s)\dd^2 F|_{u(s)}\big(v_i(s),v_j(s)\big)\nonumber\\
	&+\sum_{i=2}^n\frac{a_1(s)a_i(s)}{\lambda_i(s)-\lambda_1(s)}z_i^*(s)\dd^2 F|_{u(s)}\big(v_1(s),v_1(s)\big)
\end{align}
where $v_i(s)$, $i=1,\ldots,n$ are defined as in \eqref{0.0.2}.

In what follows, we estimate each term that appears in \eqref{4.1}. 

By the $C^2$-continuity of $\gamma$, $\big\langle\ddot{\gamma}(s),z_1(s)\big\rangle$ is uniformly bounded over $s\in I$.

By definition of $f(s)$ in \eqref{eq2:2} we have
\begin{align*}
    &\frac{a_1(s)}{\sqrt{\lambda_1(s)}}\left(\sum_{i=2}^n\frac{a_i(s)\sqrt{\lambda_i(s)}}{\lambda_i(s)-\lambda_1(s)}z_1^*(s)\dd^2 F|_{u(s)}\big(v_i(s),v_1(s)\big)\right)\\
    =&\frac{a_1(s)}{\sqrt{\lambda_1(s)}}\left(f(s)+\sum_{i=2}^n\frac{a_i(s)}{\sqrt{\lambda_i(s)}}\frac{\lambda_1(s)}{\lambda_i(s)-\lambda_1(s)}z_1^*(s)\dd^2 F|_{u(s)}\big(v_i(s),v_1(s)\big)\right)\\
    =&\frac{a_1(s)}{\sqrt{\lambda_1(s)}}f(s)+\sum_{i=2}^n\sqrt{\frac{\lambda_1(s)}{\lambda_i(s)}}\frac{a_1(s)a_i(s)}{\lambda_i(s)-\lambda_1(s)}z_1^*(s)\dd^2 F|_{u(s)}\big(v_i(s),v_1(s)\big)
\end{align*}

Since $\lambda_1$ is bounded over $I$ (see Proposition~\ref{p3.1}), using Assumption (A), we have
\begin{align}\label{4.1.1}
    \left|\sum_{i=2}^n\sqrt{\frac{\lambda_1(s)}{\lambda_i(s)}}\frac{a_1(s)a_i(s)}{\lambda_i(s)-\lambda_1(s)}z_1^*(s)\dd^2 F|_{u(s)}\big(v_i(s),v_1(s)\big)\right|\leqs &\sup_{i=1,\cdots,n;~s\in I}a_i(s)^2\cdot \frac{(n-1)C}{\lambda_i(s)-\lambda_1(s)}\nonumber\\
    \leqs&\sup_{i=1,\cdots,n;~s\in I}a_i(s)^2\cdot \frac{2(n-1)C}{\lambda_0}
\end{align}
for $s\in I$, where $C$ is the constant of condition \eqref{eq:1}. Similarly,
\begin{align*}
    &\left|\sum_{i,j=2}^n\frac{a_i(s)a_j(s)\sqrt{\lambda_1(s)}}{\big(\lambda_i(s)-\lambda_1(s)\big)\sqrt{\lambda_j(s)}}z_i^*(s)\dd^2 F|_{u(s)}\big(v_i(s),v_j(s)\big)\right|\leqs\sup_{i=2,\cdots,n;~s\in I}a_i(s)^2\frac{2(n-1)^2C}{\lambda_0},\\
    &\left|\sum_{i,j=2}^n\frac{a_i(s)a_j(s)\sqrt{\lambda_i(s)}}{\big(\lambda_i(s)-\lambda_1(s)\big)\sqrt{\lambda_j(s)}}z_1^*(s)\dd^2F|_{u(s)}\big(v_i(s),v_j(s)\big)\right|\leqs\sup_{i=2,\cdots,n;~s\in I}a_i(s)^2\frac{2(n-1)^2C}{\lambda_0},
\end{align*}
for $s\in I$, while for the last term in \eqref{4.1} we have the same bound as in \eqref{4.1.1}. The conclusion follows.
\hfill$\Box$

\section{Proof of the Main Theorems}\label{S3}
According to the discussion in the introduction, it is enough to prove that every $C^2$ continuous path $\gamma$ on $\mathbb{R}^n$ can be lifted by $F$ (after Definition~\ref{d1.1}). For that, it suffices to prove that all the PLEs \eqref{1.1}
have global solutions on $[0,1]$. Finally we only need to prove that such solutions remain bounded on their intervals of existence with a bound independent of $s$.
\subsection{Proof of Theorem~\ref{t1}}
If $\|u(s)\|_X\leqs R$ for all $s\in I$, then the PLE has a global solution since its right-hand side is actually bounded over its interval of definition.
Otherwise, choose $s_0$, $s_1$ such that  $\bar{s}\leqs s_0<s_1\leqs\sigma_0$ and $\|u(s)\|_X> R$ for all $s\in(s_0,s_1)$. 

We shall show that under the assumption of the proposition, the difference between $\|u(s_0)\|_X$ and $\|u(s_1)\|_X$ is uniformly bounded. As $a_1(\cdot)$ is bounded over $I$,  \eqref{3.5.0} implies that there exist $C_1$, $C_2\geqs0$ such that
	 \begin{align}\label{4.2.1}
     \begin{array}{rcl}
          \left|\int_{s_0}^{s_1}a_1(s)\frac{\dd\sqrt{\lambda_1}}{\dd s}(s)\dd s\right|\!\!\!\!&\geqs&\!\!\!\!\left|\int_{s_0}^{s_1}\frac{a_1(s)^2}{\sqrt{\lambda_1(s)}}h(s)\dd s\right|-C_1\\
	 	\sqrt{\lambda_1(s)}\Big|\frac{\dd a_1}{\dd s}(s)\Big|\!\!\!\!&\leqs& \!\!\!\!C_2
     \end{array}
	 \end{align}
	 By the boundedness of $a_1(s)\sqrt{\lambda_1(s)}$ on $I$,  we have from \eqref{4.2} that
	 \begin{align}\label{4.4}
	 	\left|\int_{s_0}^{s_1}\frac{a_1(s)^2}{\sqrt{\lambda_1(s)}}h(s)\dd s\right|\leqs&\left|
        \int_{s_0}^{s_1}a_1(s)\frac{\dd\sqrt{\lambda_1}}{\dd s}(s)\dd s\right|+C_1\nonumber\\
        \leqs&\left|\Big[a_1\sqrt{\lambda_1}\Big]_{s=s_0}^{s=s_1}-\int_{s_0}^{s_1}\sqrt{\lambda_1(s)}\frac{\dd a_1}{\dd s}(s)\dd s\right|+C_1\nonumber\\
        \leqs&\Big|a_1(s_1)\sqrt{\lambda_1(s_1)}-a_1(s_0)\sqrt{\lambda_1(s_0)}\Big|+C_2(s_1-s_0)+C_1\leqs C_3,
	 \end{align}
	 for some positive constant $C_3$.

     By assumption \eqref{eq:parent}, $|h(s)|\leqs C$ for all $s\in I$, hence combining \eqref{4.4} and \eqref{3.00} with the upper bounds of $\lambda_1(s)$ and $|a_1(s)|$, we see that there exists $C_4\geqs0$ such that
	 \begin{align}\label{4.5}
	 	\int_{s_0}^{s_1}\Big\|\frac{\partial u(s)}{\partial s}\Big\|_X^2\sqrt{\lambda_1(s)}|h(s)|\dd s\leqs&\int_{s_0}^{s_1}\Big(\frac{|a_1(s)|}{\sqrt{\lambda_1(s)}}+\bar{C}\Big)^2\sqrt{\lambda_1(s)}|h(s)|\dd s\nonumber\\
        \leqs&\int_{s_0}^{s_1}\frac{a_1(s)^2}{\sqrt{\lambda_1(s)}}|h(s)|+2\bar{C}|a_1(s)h(s)|+\bar{C}^2\sqrt{\lambda_1(s)}|h(s)|\dd s\nonumber\\
        \leqs&C_3+2\bar{C}C(s_1-s_0)\sup\limits_{s\in I}|a_1(s)|+\frac{\sqrt{\lambda_0}}{2}\bar{C^2}C(s_1-s_0)
        \leqs C_4,
	 \end{align} 
	 Finally, since by definition \eqref{eq2:1}, and by using  assumption \eqref{eq:2}, \eqref{4.5} yields that  
    \begin{align}\label{4.7}
        \int_{s_0}^{s_1}\Big\|\frac{\partial u(s)}{\partial s}\Big\|_X^2\frac{1}{\xi(\|u(s)\|_X)^2}\dd s\leqs \int_{s_0}^{s_1}\Big\|\frac{\partial u(s)}{\partial s}\Big\|_X^2\sqrt{\lambda_1(s)}|h(s)|\dd s\leqs C_4
    \end{align}
    which implies that 
    \begin{align}\label{4.8}
        H(\|u(s_0)\|)-H(\|u(s_1)\|)=\int_{\|u(s_0)\|}^{\|u(s_1)\|}\frac{1}{\xi(\|u(s)\|)}\dd \|u(s)\|_X\leqs\int_{s_0}^{s_1}\Big\|\frac{\partial u(s)}{\partial s}\Big\|\frac{1}{\xi(\|u(s)\|)}\dd s\leqs C_5
    \end{align}
     where $H:[R,+\infty)\ra\R^+$ is defined as $H(x)=\int_R^{x}\frac{\dd t}{\xi(t)}$ and the second last inequality is obtained by the Cauchy-Schwarz  inequality and \eqref{4.7}. By assumption on $\xi$, $H$ defines an homeomorphism from $[R,\infty)$ to $\mathbb{R}^+$; we denote its inverse by $H^{-1}$. If $[s_0,s_1]=I$, then $\|u(s))\|\leq H^{-1}(C_5+H(\|u^0)\|))$  for $s\in I$ and the conclusion follows.
     If not, then it holds that 
     $\|u(s_0)\|=R$ (or $\|u(s_1)\|=R$) 
     and hence $\|u(s))\|\leq H^{-1}(C_5+H(R))$  for $s\in (s_0,s_1)$. 
     Therefore, $u$ is again bounded over $I$ with a bound independent of $s\in I$. The proof of Theorem~\ref{t1} is complete.
    \hfill$\Box$

    \begin{remark}
        In the special case as described in Remark \ref{r2.0.0}, the estimate \eqref{4.7} and the condition \eqref{0.0.01} imply that there exists $C_4'\geqs0$ such that
	 \begin{align}\label{4.6}
	 	\int_{s_0}^{s_1}\Big\|\frac{\partial u(s)}{\partial s}\Big\|_X^2\frac{K_1}{K_2\|u(s)\|_X^2}\dd s\leqs \int_{s_0}^{s_1}\Big\|\frac{\partial u(s)}{\partial s}\Big\|_X^2\sqrt{\lambda_1(s)}|h(s)|\dd s\leqs C_4'
	 \end{align}
	 with $s_0$, $s_1$ chosen as above, which implies that there exists $C_5'\geqs0$ such that
	 \begin{align*}
	 	\int_{s_0}^{s_1}\Big\|\frac{\partial u(s)}{\partial s}\Big\|_X\frac{1}{2\|u(s)\|_X}\dd s\leqs C_5'.
	 \end{align*}
	 Hence calculating the variation of $\ln\|u(s)\|_X$ on the interval yields
	 \begin{align}\label{4.0.4}
	 	\ln\|u(s_1)\|_X-\ln\|u(s_0)\|_X=&\int_{s_0}^{s_1}\frac{\dd}{\dd s}\ln(\|u(s)\|_X)\dd s\nonumber\\
	 	=&\int_{s_0}^{s_1}\Big\langle\frac{\partial u(s)}{\partial s}, \frac{u(s)}{2\|u(s)\|^2}\Big\rangle_X\dd s\nonumber\\
	 	\leqs&\int_{s_0}^{s_1}\Big\|\frac{\partial u(s)}{\partial s}\Big\|_X\frac{1}{2\|u(s)\|_X}\dd s\nonumber\\
        \leqs& C_5'
	 \end{align}
	 which proves the proposition by showing the finite variation of $\|u(s)\|$ on $s\in[\bar{s},\sigma)$. Moreover, under the condition of Theorem~\ref{t1}, by Lemma \ref{l3.2} and the continuity of $a_1(s)$ and $h(s)$, we have
        \begin{align*}
			\Big\|\frac{\partial u(s)}{\partial s}\Big\|_X\leqs \frac{|a_1(s)|}{\sqrt{\lambda_1(s)}}+\bar{C}\leqs\frac{1}{|h(s)|\sqrt{\lambda_1(s)}}|a_1(s)h(s)|+\bar{C}\leqs\frac{K_2\|u(s)\|^2_X}{K_1}|a_1(s)h(s)|+\bar{C}\leqs \bar{C}_1,
		\end{align*}
        where $\bar{C}_1\geqs0$ is a constant, implying that {there is no blow-up on the right hand side of the PLE.}
	\end{remark}


    

\subsection{Proof of Theorem~\ref{t1.1}}
    We next prove the boundedness of solutions of the PLE as they hit 
    $\tilde{S}$ at $s=1$. With no loss of generality, we assume that Assumption (A') holds (otherwise we are done). 

Similarly as in the proof of Theorem~\ref{t1}, it suffices to consider the existence of the solution on $s\in I$ when $\|u(s)\|_X> R$. Choose $s_0$, $s_1$ such that  $s_0<s_1$ in $I$ and $\|u(s)\|_X>R$ for $s\in (s_0,s_1)$.


Set $g(s)=\frac{a_1(s)}{\sqrt{\lambda_1(s)}}$ for $s\in I$. Recalling the estimates from Proposition \ref{p3.1} and \ref{p3.2}, we have 
\begin{align}\label{4.0.6}
    \frac{\dd g}{\dd s}(s)&=-\frac{1}{\lambda_1(s)}\left(\sqrt{\lambda_1(s)}\frac{\dd a_1}{\dd s}(s)-a_1(s)\frac{\dd\sqrt{\lambda_1}}{\dd s}(s)\right)\nonumber\\
    &=-\frac{1}{\lambda_1(s)}\left(a_1(s)f(s)+O\big(\sqrt{\lambda_1(s)}\big)-\frac{a_1(s)^2}{\sqrt{\lambda_1(s)}}h(s)-a_1(s)f(s)\right)\nonumber\\
    &=\frac{1}{\sqrt{\lambda_1(s)}}\Big(g(s)^2h(s)+O(1)\Big),
\end{align}
 implying that, after using condition \eqref{eq1:2}, there exists a positive constant $\tilde{C}$ such that 
\begin{align}\label{4.0.5}
    \Big|\frac{\dd g}{\dd s}(s)\Big|\geqs\frac{1}{\sqrt{\lambda_1(s)}}\Big(Kg(s)^2-\tilde{C}\Big).
\end{align}

Consider now the open subset $E_+$ of $(s_0,s_1)$ defined as
\[
E_+=\Big\{s\in (s_0,s_1)~\Big|~ |g(s)|>\sqrt{\frac{2\tilde{C}}{K}}\Big\}.
\]
Then $E_+$ is the disjoint union of at most countable many intervals $I_j=(a_j,b_j)$, $j\in J\subset \mathbb{N}$. On each $I_j$, we have
\begin{align}\label{4.0.9}
    \Big|\frac{\dd g}{\dd s}(s)\Big|>\frac{K}{2\sqrt{\lambda_1(s)}}g(s)^2>0
\end{align}
which implies that not only $g$ has a constant sign on $I_j$ but $\frac{\dd g}{\dd s}(s)$ as well, implying that
$g$ is strictly monotone and hence either 
$|g(a_i)|$ or $|g(b_j)|$ is larger than 
$\sqrt{\frac{2\tilde{C}}{K}}$. Therefore
for $j\in J$, either $a_i=s_0$ or $b_j=s_1$. One deduces that either $J$ is empty, i.e., 
$|g(s)|\leqs \sqrt{\frac{2\tilde{C}}{K}}$
for $s\in (s_0,s_1)$ or $J$ has at most two elements. On such an interval (let say $(s_0,b_1)$),   using \eqref{4.0.9}, one has
\begin{align}
\int_{s_0}^{b_1}|g(s)|\dd s\leqs
    C_0\int_{s_0}^{b_1}\frac{\dd s}{\sqrt{\lambda_1(s)}}\leqs \frac{2C_0}{K}\Big|\int_{s_0}^{b_1}\frac{\dd g}{\dd s}(s)\frac{\dd s}{g(s)^2}\Big|\leqs \frac{2C_0}{K}\Big|\Big(\frac{1}{g(s_0)}-\frac{1}{g(b_1)}\Big)\Big|\leqs \frac{2\sqrt{2}C_0}{\sqrt{\tilde{C}K}}
\end{align}
One deduces at once that $g$ is integrable over $(s_0,s_1)$. Therefore, by Lemma \ref{l3.2} we see that the total variation of $\|u(s)\|$ is bounded on $[s_0,s_1]$. Summarizing the above arguments, we have the uniform boundedness of $\|u(s)\|_X$ for $s\in[\bar{s},\sigma_0)$. The proof of Theorem~\ref{t1.1} is complete.
\hfill$\Box$

\section{Application and Conclusion}\label{S4}

    An important application of Theorem \ref{t1} is that it justified the well-posedness of the homotopy continuation method for the motion planning of nonlinear control systems over $\R^n$ \cite{chs}. Consider the following nonlinear control system
    \begin{align}\label{5.0.1}
        \dot{x}(t)=f(x(t),u(t))
    \end{align}
    where $f(\cdot,u)$ is a $C^3$ vector field on $\R^n$ for all $u\in\R^m$, and is $C^2$ with respect to $u$.
    For any $\bar{x}\in \R^n$, $\bar{u}\in L^2([0,T],\R^m)$, denote by $\varphi_{\bar{u}}(\bar{x},t)$ the solution of  (\ref{5.0.1}) with $u=\bar{u}$ and initial value $\bar{x}$, $t\in[0,T]$. Let $x$, $x_1$ be two points on $\R^n$. The motion planning problem aims at finding a control $u^*\in L^2([0,T],\R^m)$ such that $\varphi_{u^*}(x,T)=x_1$. 
	
	For $T>0$, define the {\it endpoint map} corresponding to (\ref{5.0.1}) starting from $x$ in time $T$ as
the mapping $\E_{x,T}$ which associates to a control $u\in L^2([0,T],\R^m)$ the terminal point $\varphi_u(x,T)$.
In case $\E_{x,T}$ well-defined for all $u\in L^2([0,T],\R^m)$, it holds that $\E_{x,T}:L^2([0,T],\R^m)\to \mathbb{R}^n$ is a $C^2$ map \cite{tre} and the  controllability of (\ref{5.0.1}) from $x_0$ in time $T$
is equivalent to the surjectivity of $\E_{x,T}$. 

    One can see that the solvability of the motion planning is equivalent to the (local) surjectivity of the endpoint map. Therefore, the Hadamard-L\'{e}vy type theorem, Lemma \ref{l1.0}, can be applied to generate a solution. If condition \eqref{0.11} holds, then solving the PLE \eqref{1.1} with the $F$ being the endpoint map $\E_{x,T}$ will give the desired control as $u(1)$. In \cite{ch96,su93}, condition \eqref{0.11} is satisfied under strong conditions on the Lie configuration of the nonholonomic dynamics, which, however, are not easy to verify. Based on these results, numerical methods were established for nonlinear motion planning \cite{alc}. In \cite{sch} and \cite{ji} the authors proposed a regularization method to avoid the singularity but did not show the well-posedness the solution at the limit of the regularization parameter.

    On the other hand, our main result, Theorem \ref{t1.1}, implies the following sufficient condition for the well-posedness of the  motion planning problem.

    \begin{corollary}\label{c4.2}
        Let $x,x_1\in\R^n$. Consider the control system \eqref{5.0.1} and its endpoint map $\E_{x,T}$. Let $u^0\in L^2([0,T],\R^m)$ and $\gamma:[0,1]\ra \R^n$ be a $C^1$ curve satisfying $\gamma(0)=\E_{x,T}(u^0)$ and $\gamma(1)=x_1\in S^x_T$, where $S^x_T$ is the set of singular values of $\E_{x,T}$. If $\gamma$ and $\E_{x,T}$ satisfy the conditions in Theorem \ref{t1.1}, then there exists an $L^2$-bounded control $u_1$ such that $\E_{x,T}(u_1)=x_1$.
    \end{corollary}


    \begin{remark}
        In the context of geometric control theory, condition \eqref{eq1:2} in our main theorem has an intuitive interpretation when paraphrased as
        \begin{align*}
            |z^*\dd^2\E_{x,T}(v,v)|\geqs K, \quad\forall z\in\R^n, ~\forall v\in L^2([0,T],\R^m), ~\|z\|_{\R^n}=1, ~\|v\|_{L^2}=1.
        \end{align*}
        If $v$ lies in the singular set of the endpoint map, then the above condition implies that $\E_{x,T}(v)$ is not a {\it conjugate point} of the system \eqref{5.0.1} \cite{bon}. Therefore \eqref{eq1:2} can be interpreted as a generalized condition for non-conjugacy.
    \end{remark}

    In conclusion, by analyzing the spectral of the Gramian along the PLE, we proved the global well-posedness of its solution, and further showed the surjectivity of the given map via its path-lifting property, of which the nonlinear motion planning problem can be viewed as a special case. These results broaden the applicability of the homotopy continuation method and lay the groundwork for further generalizations.
    
    For future work, we propose the following directions:
    \begin{enumerate}
        \item Extension to infinite-dimensional target spaces. This would involve estimating  the differential of the spectral of the Gramian as depicted in Section \ref{S2} over infinite-dimensional manifolds.
        \item Relaxation of assumption (A) on the lower bound of the eigenvalues of the Gramian. Eliminating this assumption on singularity will improve generality and requires differential analysis of the canonical projectors of the Gramian.
        \item Higher-order generalizations. One sees that the analysis in Section \ref{S2} relies on estimating the first order differentiation of the spectral; when the first derivative of the eigenvalue vanishes, higher-order differentials may govern the behavior of the PLE, opening possibilities to further refinement of surjectivity conditions. 
        
    \end{enumerate}


	
%

\end{document}